\documentclass[siam10.clo]{siamltex}

\begin{document}
\title{A different view on third and fourth order efficiency}

\author{ Shanti Venetiaan\thanks{Institute for Graduate Studies and Research, Anton de Kom Universiteit van Suriname, Paramaribo, Suriname. E-mail: s.venetiaan@uvs.edu} }

\maketitle

\begin{abstract}
In this paper a different approach to consider third and fourth order efficiency is suggested for non symmetric cases. Fourth order efficiency does not follow automatically but only after some adjustments are made.
\end{abstract}



\pagestyle{myheadings} \thispagestyle{plain} \markboth{ SHANTI VENETIAAN}{{\small THIRD AND FOURTH ORDER EFFICIENCY}}

\section{Introduction}
The famous phenomenon "first order efficiency implies second order efficiency", discovered by Pfanzagl(1979) and later studied by others like Bickel, Chibisov and van Zwet(1981), Klaassen and Venetiaan (1994) has led to the study of third order efficiency implying fourth order efficiency, a conjecture of Ghosh(1994). Akahira(1996) did so by studying the concentration probability of symmetric intervals. Klaassen and Venetiaan (2009) look at this issue differently by considering the so-called confidence interval inequality and show that the distribution function of the maximum likelihood estimator (MLE) attains this bound up to third order and as a consequence also up to fourth order, because of the symmetry in the expansions. In this paper the author goes back to the confidence interval inequality and the distribution function of the MLE. The same framework is used as in Klaassen and Venetiaan (2009). This means that we consider estimation of the location parameter $\theta$ in the one-dimensional location model of independently and identically distributed random variables $X_1,...,X_n$. We assume that the common distribution of the random variables has finite Fisher information for location, and as a consequence this distribution is absolutely continuous with an absolutely continuous density $f(\cdot -\theta)$ with derivative $f'(\cdot -\theta)$. Without loss of generality we may assume that $I(f) = \int (f'/f)^2 f=1$. Furthermore we name the distribution function of the MLE($T_n$), $G_n$, namely
\begin{equation}
G_n(y) = P_0(\sqrt{nI(f)} T_n \le y), \quad y\in R
\end{equation}

In Klaassen and Venetiaan(2009) an expansion for the normalized length $G_n^{-1}(1-\alpha/2)-G_n^{-1}(\alpha/2)$ of the symmetric confidence interval based on the maximum likelihood estimator is derived. In the same paper a lower bound $\epsilon$ is introduced for confidence intervals and it is shown that for symmetric intervals the normalized length $G_n^{-1}(1-\alpha/2)-G_n^{-1}(\alpha/2)$ and $\epsilon$ coincide up to the $1/n$-term and because of symmetry of the polynomials in the expansion and symmetry of the interval, the $1/n\sqrt n$-term. This enables us to say that third order efficiency is obtained and automatically fourth order efficiency. The fourth order efficiency result is entirely based on the fact that symmetric confidence intervals are studied. In the present paper the author considers one side of the confidence intervals, to avoid the symmetry and hence the automatic efficiency in fourth order.  The expansion for the bound in the confidence interval inequality, $\epsilon$ is split into two sides $\epsilon_v$ and $\epsilon_u$, which both have the same kind of terms . Then we put $v=1-\alpha/2$ and $u=\alpha/2$, study $G_n^{-1}(1-\alpha/2)-G_n^{-1}(1/2)-\epsilon_{1-\alpha/2}$ and find that the third order term vanishes, but not the fourth order term. 
 In Section 2 we will state the relevant expansions and perform the splitting of $\epsilon$. We also give the result of the comparison. In Section 3 the author suggests an adjustment which may result in fourth order efficiency. The reader will notice that this paper does not give any formal proofs and there is no logical explanation for splitting up $\epsilon$, but the author just wanted to mention the peculiarities.

\section{Confidence intervals, one side}

We will be using the following notation.
\begin{eqnarray} \eta_2 &= & E\psi_2^2(X_1), \eta_3 = E\psi_1^3(X_1), \eta_4=E\psi_1^4(X_1),\nonumber\\
\eta_5&=& E\psi_1^5(X_1), \eta_6=E(\psi_2(X_1)\psi_3(X_1))\nonumber\\
 \quad& &\mbox{with} \quad\psi_i(x) = \frac{f^{(i)}}{f}(x).\nonumber\end{eqnarray}

Recall that in Venetiaan(2009) an expansion was derived for $G_n^{-1}(v)$. Note that 
\begin{eqnarray}G_n^{-1}(v)-G_n^{-1}(1/2)&=&\Phi^{-1}(v)\nonumber\\
&\quad+&\frac{\eta_3}{12\sqrt n}(\Phi^{-1}(v))^2\nonumber\\
&\quad+&\frac1{n}\biggl[(-\frac{5\eta_4}{72}-\frac18+\frac{\eta_2}{6}-\frac{\eta_3^2}{72})(\Phi^{-1}(v)^3+ (-\frac{\eta_3^2}{36}-\frac18+\frac{\eta_4}{24})\Phi^{-1}(v)\biggr]\nonumber\\
&\quad+& \frac1{n\sqrt n}\biggl[(\frac{\eta_2\eta_3}{24}-\frac{\eta_3\eta_4}{144}-\frac{\eta_3}{48}-\frac{\eta_6}{8}+\frac{\eta_5}{30}-\frac{19\eta_3^3}{1728})((\Phi^{-1}(v))^4\nonumber\\
&\qquad+&(\frac{\eta_3\eta_4}{48}-\frac{67\eta_3^3}{1296}-\frac{5\eta_3}{48}+\frac{\eta_2\eta_3}{12}-\frac{\eta_5}{80})(\Phi^{-1}(v))^2\biggr]\nonumber\\
&\quad+& o(\frac1{n\sqrt n})\label{eqno5}
\end{eqnarray}
 The subtraction of $G_n^{-1}(1/2)$ is done to get the expansion for the median unbiased case.

On the other hand Klaassen and Venetiaan (2009) introduced a bound for confidence intervals and derived expansions for this bound, namely
\begin{eqnarray}
\epsilon &=& \Phi^{-1}(v)-\Phi^{-1}(u)\nonumber\\
&+&\frac{\eta_3}{12\sqrt n}((\Phi^{-1}(v)^2-(\Phi^{-1}(u))^2)\nonumber\\
&+& \frac1{n}\biggl[(\frac{\eta_2}{24}+\frac{5\eta_3^2}{288}-\frac{\eta_4}{36})((\Phi^{-1}(v))^3-(\Phi^{-1}(u))^3)\nonumber\\
&+&(\frac18-\frac{\eta_2}{8}+\frac{\eta_3^2}{32}+\frac{\eta_4}{24})(\Phi^{-1}(v)^2\Phi^{-1}(u)-\Phi^{-1}(u)^2\Phi^{-1}(v))\nonumber\\
&+& (\frac{\eta_4}{24}-\frac18-\frac{\eta_3^2}{36})(\Phi^{-1}(v)-\Phi^{-1}(u))\biggr]\nonumber\\
&+& \frac1{n\sqrt n}\biggl[(\frac{\eta_2\eta_3}{48}+\frac{\eta_3^3}{216}-\frac{\eta_3\eta_4}{72}+\frac{\eta_5}{80}-\frac{\eta_6}{24})(\Phi^{-1}(v)^4-\Phi^{-1}(u)^4)\nonumber\\
&+&(\frac{\eta_3}{24}-\frac{\eta_2\eta_3}{24}+\frac{\eta_3^3}{48}-\frac{\eta_5}{48}+\frac{\eta_6}{12})(\Phi^{-1}(v)^3\Phi^{-1}(u)-\Phi^{-1}(u)^3\Phi^{-1}(v))\nonumber\\
&+&(\frac{7\eta_3\eta_4}{144}-\frac{\eta_3}{48}-\frac{5\eta_3^3}{162}-\frac{\eta_5}{80})(\Phi^{-1}(v)^2-\Phi^{-1}(u)^2)\biggr]+o(\frac1{n\sqrt n})
\end{eqnarray}

Now , when we consider $G_n^{-1}(1-\alpha+u)-G_n^{-1}(u)-\epsilon$, the $1/\sqrt n$- and $1/n\sqrt n$-term will automatically vanish when we choose the $u$ which optimizes in first order, namely $u= \alpha/2$. This is because $\Phi^{-1}(1-\alpha/2) = - \Phi^{-1}(\alpha/2)$. So it seems that second and fourth order efficiency are automatically obtained when the symmetric case is studied, because the corresponding terms have even polynomials. But, as we know the phenomenon "first order efficiency implies second order efficiency" also works in non symmetric cases. We see that even the choice of $u$ is not of influence to obtain that result. 
A close look at $\epsilon$ shows that there is a symmetric structure in there and we will split $\epsilon$ in $\epsilon_v$ and $\epsilon_u$. As a consequence 
\begin{eqnarray}
G_n^{-1}(v) &-& G_n^{-1}(u) - \epsilon =\nonumber\\
G_n^{-1}(v)&-& G_n^{-1}(1/2) - (G_n^{-1}(u)-G_n^{-1}(1/2)) - (\epsilon_v -\epsilon_u)=\nonumber\\
(G_n^{-1}(v) &-& G_n^{-1}(1/2)-\epsilon_v) -(G_n^{-1}(u)-G_n^{-1}(1/2)-\epsilon_u)\label{eqno4}
\end{eqnarray}
Now from here on we will just study one side of the expression in (\ref{eqno4})
with 
\begin{eqnarray}
\epsilon_v &=& \Phi^{-1}(v)\nonumber\\
&+&\frac{\eta_3}{12\sqrt n}((\Phi^{-1}(v)^2)\nonumber\\
&+& \frac1{n}\biggl[(\frac{\eta_2}{24}+\frac{5\eta_3^2}{288}-\frac{\eta_4}{36})((\Phi^{-1}(v))^3)\nonumber\\
&+&(\frac18-\frac{\eta_2}{8}+\frac{\eta_3^2}{32}+\frac{\eta_4}{24})(\Phi^{-1}(v)^2\Phi^{-1}(u))]\nonumber\\
&+& (\frac{\eta_4}{24}-\frac18-\frac{\eta_3^2}{36})(\Phi^{-1}(v)\biggr]\nonumber\\
&+& \frac1{n\sqrt n}\biggl[(\frac{\eta_2\eta_3}{48}+\frac{\eta_3^3}{216}-\frac{\eta_3\eta_4}{72}+\frac{\eta_5}{80}-\frac{\eta_6}{24})(\Phi^{-1}(v)^4)\nonumber\\
&+&(\frac{\eta_3}{24}-\frac{\eta_2\eta_3}{24}+\frac{\eta_3^3}{48}-\frac{\eta_5}{48}+\frac{\eta_6}{12})(\Phi^{-1}(v)^3\Phi^{-1}(u))\nonumber\\
&+&(\frac{7\eta_3\eta_4}{144}-\frac{\eta_3}{48}-\frac{5\eta_3^3}{162}-\frac{\eta_5}{80})(\Phi^{-1}(v)^2)\biggr]\label{eqno3}
\end{eqnarray}

We substitute $v = 1-\alpha/2$ and $u =\alpha/2$ and put $z=\Phi^{-1}(1-\alpha/2)$, then (\ref{eqno3}) becomes
\begin{eqnarray}
\epsilon_v &=&z +\frac{\eta_3}{12\sqrt n}z^2\nonumber\\
&+& \frac1{n}\biggl[(-\frac{5\eta_4}{72}-\frac18+\frac{\eta_2}{6}-\frac{\eta_3^2}{72})z^3 + (\frac{\eta_4}{24}-\frac18-\frac{\eta_3^2}{36})z\biggr]\nonumber\\
&+& \frac1{n\sqrt n}\biggl[(\frac{\eta_2\eta_3}{16}-\frac{7\eta_3^3}{432}-\frac{\eta_3\eta_4}{72}+\frac{\eta_5}{30}-\frac{\eta_6}{8})z^4
+(\frac{7\eta_3\eta_4}{144}-\frac{\eta_3}{48}-\frac{5\eta_3^3}{162}-\frac{\eta_5}{80})z^2\biggr]+o(\frac1{n\sqrt n})
\end{eqnarray}

We also plug in $1-\alpha/2$ for $v$ in (\ref{eqno5}) and then the left side of the expression in(\ref{eqno4}) becomes
\begin{eqnarray}
&&\frac1{n\sqrt n}(\frac{\eta_3}{48}-\frac{\eta_2\eta_3}{48}+\frac{\eta_3\eta_4}{144}+\frac{\eta_3^3}{192})z^4+(\frac{\eta_2\eta_3}{12}-\frac{\eta_3}{12}-\frac{\eta_3\eta_4}{36}-\frac{\eta_3^3}{48})z^2)+o(\frac1{n\sqrt n})\\
&=& \frac{\eta_3}{48n\sqrt n}(1-\eta_2+\frac{\eta_4}{3}+\frac{\eta_3^2}{4})(z^4-4z^2)+o(\frac1{n\sqrt n}),
\end{eqnarray}

namely, the $\frac{1}{n}$-term vanished, in other words, third order efficiency is obtained. 

It seems that in this framework, fourth order efficiency for the non symmetric case may only be obtained if 
\begin{equation}
1-\eta_2 +\frac{\eta_4}{3}+\frac{\eta_3^2}{4}=0
\end{equation}
Let's name $1-\eta_2 +\frac{\eta_4}{3}+\frac{\eta_3^2}{4}=W$.
Via the Cauchy-Schwartz inequality it may be verified that $W\le 0$, by noting that
\begin{equation}
E\psi_1^2\psi_2(X_1)=\frac23E\psi_1^4(X_1)
\end{equation}
 and
\begin{equation}
(E\psi_1^3(X_1))^2= (-2E\{\psi_1(X_1)(\psi_1'(X_1)-E\psi_1'(X_1)) \})^2 \le 4E\psi_1^2(X_1)var(\psi_1'(X_1)).
\end{equation}

\subsection{Another view}
We study $G_n^{-1}(v)-G_n^{-1}(u) - \epsilon$ and find
\begin{eqnarray}
G_n^{-1}(v) -G_n^{-1}(u) -\epsilon &=& -\frac{W}{8n}((\Phi^{-1}(v))^3-(\Phi^{-1}(u))^3+(\Phi^{-1}(v))^2\Phi^{-1}(u)-(\Phi^{_1}(u))^2\Phi^{-1}(v))\nonumber\\
&+&\frac1{n\sqrt n}[(-\frac{\eta_3W}{48}-\frac{\eta_6}{12}+\frac{\eta_5}{48}-\frac{\eta_3^3}{96}+\frac{\eta_3\eta_4}{72})((\Phi^{-1}(v))^4-(\Phi^{-1}(u))^4)\nonumber\\
&+&(-\frac{\eta_3W}{24}-\frac{\eta_6}{12}+\frac{\eta_5}{48}-\frac{\eta_3^3}{96}+\frac{\eta_3\eta_4}{72})((\Phi^{-1}(v))^3\Phi^{-1}(u)-(\Phi^{-1}(u))^3\Phi^{-1}(v))\nonumber\\
&+&(-\frac{\eta_3W}{12})((\Phi^{-1}(v))^2-(\Phi^{-1}(u))^2)]+o(\frac1{n\sqrt n})
\end{eqnarray}
The third order term will vanish if $W=0$,i.e. third order efficiency will be obtained when $W=0$ for whatever choice of $u$ and $v$ and note that the first and second term of the $1/n\sqrt n$-term will have the same coefficient and the last one will disappear.

\section{Adding "missing" terms}
In deriving $\epsilon$ we obtain a bound for a two-sided confidence interval. It could be possible that if we were able to study one-sided intervals, we might see some extra terms appear which vanish when the two-sided interval is considered. Possibly, we might view $0$ as
 $0= \Phi^{-1}(u)^2\Phi^{-1}(v)^2-\Phi^{-1}(u)^2\Phi^{-1}(v)^2$ and $0 = \Phi^{-1}(u)\Phi^{-1}(v) -\Phi^{-1}(u)\Phi^{-1}(v)$. This means that $\epsilon_v$ might be
\begin{eqnarray}
\epsilon_v &=& \Phi^{-1}(v)\nonumber\\
&+&\frac{\eta_3}{12\sqrt n}((\Phi^{-1}(v))^2\nonumber\\
&+& \frac1{n}[(\frac{\eta_2}{24}+\frac{5\eta_3^2}{288}-\frac{\eta_4}{36})((\Phi^{-1}(v))^3)\nonumber\\
&+&(\frac18-\frac{\eta_2}{8}+\frac{\eta_3^2}{32}+\frac{\eta_4}{24})(\Phi^{-1}(v)^2\Phi^{-1}(u))]\nonumber\\
&+& (\frac{\eta_4}{24}-\frac18-\frac{\eta_3^2}{36})(\Phi^{-1}(v))\nonumber\\
&+& \frac1{n\sqrt n}\biggl[(\frac{\eta_2\eta_3}{48}+\frac{\eta_3^3}{216}-\frac{\eta_3\eta_4}{72}+\frac{\eta_5}{80}-\frac{\eta_6}{24})(\Phi^{-1}(v)^4)\nonumber\\
&+&(\frac{\eta_3}{24}-\frac{\eta_2\eta_3}{24}+\frac{\eta_3^3}{48}-\frac{\eta_5}{48}+\frac{\eta_6}{12})(\Phi^{-1}(v)^3\Phi^{-1}(u))\nonumber\\
&+&(\frac{7\eta_3\eta_4}{144}-\frac{\eta_3}{48}-\frac{5\eta_3^3}{162}-\frac{\eta_5}{80})(\Phi^{-1}(v)^2)\nonumber\\
&-& \frac{\eta_3W}{48}(\Phi^{-1}(v))^2(\Phi^{-1}(u))^2-\frac{\eta_3W}{12}\Phi^{-1}(v)\Phi^{-1}(u)\biggr]+o(\frac1{n\sqrt n})
\end{eqnarray}

If $u=\alpha/2$ and $v = 1-\alpha/2$, $G_n^{-1}(1-\alpha/2)- G_n^{-1}(1/2)-\epsilon_{1-\alpha/2}$ will become $o(1/n\sqrt n)$, which we see as third order efficiency and fourth order efficiency being obtained.


\begin{thebibliography}{10}

\bibitem{akahira}{\sc M. Akahira}, {\em Third order efficiency implies fourth order efficiency : a resolution of the conjecture of J.K. Ghosh}, Ann. Inst. Statist. Math., vol 48, No. 2 (1996), 365 - 380. 
\bibitem{bickchizwet}{\sc P.J. Bickel, D.M. Chibisov and W.R. van Zwet}, {\em On Efficiency of first and second order}, Internat. Statist. Rev. vol 49 (1981), 169 - 175.
\bibitem{ghosh}{\sc J.K. Ghosh}, {\em Higher order asymptotics}, NSF-CBMS Regional Conference Series Probability and Statistics, vol 4 (1994), Institute of Mathematical Statistics, Hayward, California.
\bibitem{klaasvene}{\sc C.A.J. Klaassen and S.A. Venetiaan}, {\em Spread inequality and efficiency of first and second order} Asymptotic Statistics ( M. Hu{$\check{s}$}kov$\acute{a}$ and P. Mandl, eds.), Proceed. Fifth Prague Symposium, Physica verlag, 1994, pp. 341 - 348.
\bibitem{klaasvene2}{\sc C.A.J. Klaassen and S. A. Venetiaan},{\em Optimizing lengths of confidence intervals}, (2009), to be published.
\bibitem{pfan2}{\sc J. Pfanzagl}, {\em First order efficiency implies second order efficiency}, Contributions to Statistics [J. Hajek Memorial Volume] ( J. Jureckov$\acute{a}$, ed.), Academia Prague, 1979, pp. 157-196.
\bibitem{vene2}{\sc S.A. Venetiaan}, {\em An expansion for the maximum likelihood estimator of location and its distribution function},(2009) to appear in Brazilian journal of Probability and Statistics.
\end{thebibliography}
\end{document}